\documentclass[11pt]{article}

\usepackage{latexsym}
\usepackage{mathrsfs}
\usepackage{amsfonts}
\usepackage{amssymb}
\usepackage{amsmath}
\usepackage{longtable}
\usepackage{graphicx}
\usepackage{pifont}
\usepackage{algorithmic}

\newtheorem{theorem}{Theorem}

\newtheorem{corollary}[theorem]{Corollary}
\newtheorem{lemma}[theorem]{Lemma}
\newtheorem{definition}{Definition}

\newtheorem{algorithm}{Algorithm}

\newcommand {\R}        {{\mathbb R}}

\newcommand {\mat}      [1] {\left[\begin{array}{#1}}
\newcommand {\rix}          {\end{array}\right]}
\newcommand {\bfx}       {{\bf  x}}

\newcommand*{\affaddr}[1]{#1} 
\newcommand*{\affmark}[1][*]{\textsuperscript{#1}}

\begin{document}

\begin{titlepage}
\title{
A Condition Analysis of the Weighted Linear Least Squares Problem
Using Dual Norms
}

\author{
Huai-An Diao\affmark[1]\thanks{hadiao@nenu.edu.cn, hadiao78@yahoo.com},
Liming Liang\affmark[1]\affmark[2]\thanks{liangliming0809@126.com},
Sanzheng Qiao\affmark[3]\thanks{qiao@cas.mcmaster.ca}\\
\affaddr{\affmark[1]School of Mathematics and Statistics,
Northeast Normal University, \\
No. 5268 Renmin Street, Chang Chun 130024,
P.R. China.}\\
\affaddr{\affmark[2]Current address: Qiqihar Experimental High School, \\
Zhonghua West Road No. 142,
Qiqihar 161000, P.R. China.}\\
\affaddr{\affmark[3]Department of Computing and Software,
McMaster University \\
Hamilton, Ontario L8S4K1 Canada.}
}

\date{}
\end{titlepage}

\maketitle

\begin{abstract}
In this paper, based on the theory of adjoint operators and dual norms,
we define condition numbers for a linear solution function of
the weighted linear least squares problem. The explicit expressions
of the normwise and componentwise condition numbers
derived in this paper can be computed at low cost when the dimension
of the linear function is low due to dual operator theory. Moreover,
we use the augmented system to perform a componentwise perturbation
analysis of the solution and residual of the weighted linear
least squares problems. We also propose two efficient condition
number estimators. Our numerical experiments demonstrate
that our condition numbers give accurate perturbation bounds and
can reveal the conditioning of individual components of the
solution. Our condition number estimators are accurate as well
as efficient.
\end{abstract}

\bigskip


\medskip


\section{Introduction} \label{secIntro}

This paper investigates the condition of the weighted linear
least squares problem using adjoint operators and dual norms.

Given $A \in \mathbb{R}^{m \times n}$, $m \ge n$, and
$\mathbf{b} \in \mathbb{R}^m$, the weighted least squares
problem (WLS)
\begin{equation} \label{eqnWLS}
\min_{\mathbf{x} \in \mathbb{R}^n} \| A \mathbf{x} - \mathbf{b} \|_W^2
= \min_{\mathbf{x}} ( A \mathbf{x} - \mathbf{b} )^{\mathrm{T}} W
( A \mathbf{x} - \mathbf{b} ) ,
\end{equation}
where $W \in \mathbb{R}^{m \times m}$ is symmetric and positive
definite, is a generalization of the standard least squares
problem (LS)
\begin{equation}\label{eq:LS}
	\min_{\mathbf{x}} \| A \mathbf{x} - \mathbf{b} \|_2^2.
\end{equation}
When $A$ is of full column rank, the problem (\ref{eqnWLS}) has
a unique solution that can be obtained by solving the normal
equations
\begin{equation} \label{eqnNormal}
A^{\mathrm{T}}WA \mathbf{x} = A^{\mathrm{T}} W \mathbf{b}.
\end{equation}
Alternatively, the solution can be obtained by solving the
augmented system:
\begin{equation} \label{eqnAugmented}
\left[ \begin{array}{cc}
W^{-1} & A \\
A^{\mathrm{T}} & 0
\end{array} \right]
\left[ \begin{array}{c}
\mathbf{d} \\ \mathbf{x}
\end{array} \right] =
\left[ \begin{array}{c}
\mathbf{b} \\ \mathbf{0}
\end{array} \right] ,
\end{equation}
where $\mathbf{d} = W(\mathbf{b} - A\mathbf{x})$ is the weighted
residual. The systems (\ref{eqnNormal}) and (\ref{eqnAugmented})
are mathematically equivalent for solving $\mathbf{x}$.
Methods for solving the weighted least squares
problem can found in \cite{LH74}.

The weighted least squares problem (\ref{eqnWLS}) can be reduced
to the standard least squares problem by the transformations
$A \leftarrow W^{1/2} A$ and
$\mathbf{b} \leftarrow W^{1/2} \mathbf{b}$, where $W^{1/2}$
is the square root of the symmetric and positive definite
$W$, that is $W^{1/2}W^{1/2} = W$. However, the weighted least
squares problem arises from applications where $W$ varies.
Thus it is efficient to solve the weighted least squares
problem rather than transforming it into the standard least
squares problem for every $W$. For example, in the interior
point method for convex quadratic programming, a convex quadratic
programming problem is transformed into
\[
\begin{array}{rl}
\displaystyle{
\min_{y \in \mathbb{R}^m}} &
\frac{1}{2} \mathbf{y}^{\mathrm{T}} H \mathbf{y} +
\mathbf{c}^{\mathrm{T}} \mathbf{y} \\
\hbox{subject to} & A^{\mathrm{T}} \mathbf{y} = \mathbf{b} \\
& \mathbf{y} \ge 0 ,
\end{array}
\]
where $H$ is symmetric and positive semi-definite and
$A$ is of full column rank. Introducing the dual variable
$\mathbf{x}$ and applying the Newton's method to the primal-dual
equation, we get
\[
\left[ \begin{array}{cc}
M & A \\
A^{\mathrm{T}} & 0
\end{array} \right] \,
\left[ \begin{array}{r}
- \Delta \mathbf{y} \\ \Delta \mathbf{x}
\end{array} \right] =
\left[ \begin{array}{c}
\mathbf{z} \\ A^{\mathrm{T}} \mathbf{y} - \mathbf{b}
\end{array} \right] ,
\]
where $M$ is the Hessian matrix and $\Delta \mathbf{y}$ and
$\Delta \mathbf{x}$ are respectively the updates for $\mathbf{y}$ and
$\mathbf{x}$ in the Newton iteration. The matrix $M$ changes
during the iteration while $A$ is fixed. Thus, in this problem,
we need to solve a sequence of weighted least squares problems
(\ref{eqnAugmented}) with variable $W$ but fixed $A$.

Another application where the weighted least squares problem arises
is the linear regression in statistics. As we know, linear least
squares model is commonly used for linear regression assuming
the response variables have the same error variance. In practice,
however, observations may not be equally reliable. In that case,
the weighted least squares model is an improvement of the
standard least squares model.

Condition number plays an important role in numerical analysis.
It is a measurement of the sensitivity of the solution of a
problem to the perturbation of its data. In 1966, Rice \cite{Rice66}
presented a general theory of condition based on the
Fr{\'e}chet derivative defined as follows.

Let $V$ and $W$ be two Banach spaces and $U$ an open subset of
$V$. Considering an operator $f:\ U \rightarrow W$, if, for
an $x \in U$, there exists a bounded linear operator
$A_x: \ V \rightarrow W$ such that
\[
\lim_{h \to 0} \frac{\| f(x+h) - f(x) - A_x(h) \|_W}
{\| h \|_V} = 0 ,
\]
then $f$ is said to be Fr{\'e}chet differentiable at $x$ and
$A_x$ is called the Fr{\'e}chet derivative of $f$ at $x$.

Let $\psi : \ \mathbb{R}^m \rightarrow \mathbb{R}^n$.
If $\psi$ is continuous and Fr{\'e}chet differentiable in
a neighbourhood of $a \in \mathbb{R}^m$, then, from Rice's
theory, the normwise condition number of $\psi$ at $a$ is
defined by
\begin{equation} \label{eqnRice}
\mathrm{cond}_{\psi} (a) = \lim_{\epsilon \to 0}
\sup_{\| \delta a \| \le \epsilon \|a\|}
\frac{\| \psi (a + \delta a) - \psi (a) \| / \| \psi (a) \|}
{\| \delta a \| / \| a \|} =
\frac{\| \psi^{\prime} (a) \| \, \|a\|}{\| \psi (a) \|} ,
\end{equation}
where $\psi^{\prime} (a)$ is the Fr{\'e}chet derivative
of $\psi$ at $a$. The condition number defined above can
be interpreted as the ratio of the relative error
in the solution to the relative error in the input data.
Clearly, the condition number defined above is norm dependent.
Moreover, it is a global measurement, which has some shortcomings.
For example, the distribution of the perturbations in data
is not represented. Also, as pointed out in \cite{Bjorck96},
for poorly scaled (imbalanced) problems, the error can
be overly estimated by the normwise condition number.
To alleviate the shortcomings, componentwise perturbation analysis
is introduced.

The componentwise error analysis of linear systems
can be found in \cite{CI95, Rump03b}. For the linear
least squares, componentwise perturbation analysis and
error bounds are given in \cite{Bjorck91, KLR98}.
In particular, for the full rank linear least squares
problem, \cite{Gratton96} and \cite{ABG07} present the condition
number when the perturbations in both the data and solution
are measured by norms. In \cite{BG09}, the perturbation
in the data is componentwise, whereas the perturbation in
the solution can be either componentwise or normwise, leading
to componentwise or mixed condition numbers.

Here is a brief review of some perturbation analyses of the
linear least squares problem (LS) and weighted linear least
squares problem (weighted LS). For the normwise perturbation analysis,
we refer the classical paper \cite{Stewart77} and references therein.
Cucker, Diao and Wei studied the mixed and componentwise
condition numbers for LS in \cite{CDW07}. The flexible normwise
condition numbers for LS was introduced in \cite{WDQ07}. In \cite{CD07},
Cucker and Diao gave explicit expressions for normwise, mixed and
componentwise condition numbers for LS under structured perturbations.
Diao and Wei proposed and derived the weighted Frobenius normwise
condition number for LS \cite{DW07}. Recently,
Diao et al. \cite{DWWQ13} studied the normwise, mixed and
componentwise condition number for LS involving Kronecker product.
For weighed LS, the  perturbation analysis can be found in
\cite{GW92,Gulliksson95,GJW02} and references therein. Wei and Wang
studied the explicit normwise condition numbers under range
conditions \cite{WW03}. Wang et al. derived the results of
the Frobenius normwise condition numbers for  weighted LS when
the coefficient matrix is of full rank \cite{WZXL09}.
In \cite{LS09}, Li and Sun derived explicit expressions of
mixed and componentwise condition numbers of the weighted LS problem.
Yang and Wang considered the flexible normwise condition numbers for
weighed LS \cite{YW16}.

In this paper, we often use the weighted generalized inverse
of $A$ defined as follows for the weighted least squares
problem (\ref{eqnWLS}). For $A \in \mathbb{R}^{m \times n}$,
let $M \in \mathbb{R}^{m \times m}$ and $N \in \mathbb{R}^{n \times n}$
be symmetric and positive definite. If there exists an
$X \in \mathbb{R}^{n \times m}$ satisfying the following
four equations
\begin{equation} \label{eqnWGI}
AXA = A, \quad XAX = X, \quad (MAX)^{\mathrm{T}} = MAX, \quad
(NXA)^{\mathrm{T}} = NXA,
\end{equation}
then $X$ is called the weighted generalized inverse of $A$
corresponding to the weight matrices $M$ and $N$ and denoted
by $A_{M,N}^{\dagger}$ \cite{WWQ04}.

When $A$ is of full column rank, the unique solution for (\ref{eqnWLS})
is $\mathbf{x} = (A^{\mathrm{T}} W A)^{-1} A^{\mathrm{T}} W \mathbf{b}$.
Setting $X = (A^{\mathrm{T}} W A)^{-1} A^{\mathrm{T}} W$,
$M = W$ and $N = I_n$, it can be verified that $X$ satisfies
the four equations in (\ref{eqnWGI}), that is, $X$ is the
weighted generalized inverse of $A$ corresponding to the weight
matrices $W$ and $I_n$. Thus, we have
$(A^{\mathrm{T}} W A)^{-1} A^{\mathrm{T}} W = A_{W,I_n}^{\dagger}$.
The unique solution for (\ref{eqnWLS}) can be given by
\[
\mathbf{x} = (A^{\mathrm{T}} W A)^{-1} A^{\mathrm{T}} W \mathbf{b} =
A_{W,I_n}^{\dagger} \mathbf{b}.
\]
This paper studies the sensitivity of the solution $\mathbf{x}$ to
the perturbation in
the data $A$ and $\mathbf{b}$ in the problem (\ref{eqnWLS}) by
applying the condition number defined in (\ref{eqnRice}). So,
corresponding to the function $\psi$ in the definition,
we define the following function mapping the data $A$ and $\mathbf{b}$
to the solution $\mathbf{x}$:
\begin{equation} \label{eqnMapSol}
g(A, \mathbf{b}) = L^{\mathrm{T}}
(A^{\mathrm{T}} W A)^{-1} A^{\mathrm{T}} W \mathbf{b} =
L^{\mathrm{T}} A_{W,I_n}^{\dagger} \mathbf{b} ,
\end{equation}
where $L$ is an $n$-by-$k$, $k \le n$, matrix introduced for
the selection of the solution components. For example, when
$L = I_n$ ($k=n$), all the $n$ components of the solution $\mathbf{x}$
are equally selected. When $L = \mathbf{e}_1$ ($k=1$), the first unit
vector in $\mathbb{R}^n$, then only the first component of the
solution is selected.

In this paper, by using adjoint operators and dual norms,
we derive the condition numbers of the weighted
least square problem, where the perturbation in the data is
componentwise whereas the perturbation in the solution
is componentwise or normwise. After a brief review of adjoint
operators and dual norms and their application to condition
number in Section~\ref{secDual}, the condition numbers for
the weighted least squares problem
are presented in Section~\ref{secWLS}. An error
analysis of the augmented system (\ref{eqnAugmented}) is
performed in Section~\ref{secAugmented}. Finally, Section~\ref{secExp}
shows our numerical experiment results.

\section{Dual Techniques} \label{secDual}

Consider a linear operator $J : \ E \rightarrow G$,
between two Euclidean spaces $E$ and $G$ with the scalar
products $\langle \cdot , \cdot \rangle_E$ and
$\langle \cdot , \cdot \rangle_G$ respectively. Denote
the corresponding norms by $\| \cdot \|_E$ and $\| \cdot \|_G$
respectively. We first define adjoint operator and dual norm.

\begin{definition} \label{defAdjoint}
The adjoint operator $J^* : \ G \rightarrow E$ of $J$
is defined by
\[
\langle \mathbf{y}, J \mathbf{x} \rangle_G =
\langle J^* \mathbf{y}, \mathbf{x} \rangle_E ,
\]
where $\mathbf{x} \in E$ and $\mathbf{y} \in G$.
\end{definition}

\begin{definition} \label{defDualNorm}
The dual norm $\| \cdot \|_{E^*}$ of the norm $\| \cdot \|_E$
is defined by
\[
\| \mathbf{x} \|_{E^*} = \max_{\mathbf{u} \neq 0}
\frac{\langle \mathbf{x} , \mathbf{u} \rangle_E}
{\| \mathbf{u} \|_E} .
\]
The dual norm $\| \cdot \|_{G^*}$ can be defined similarly.
\end{definition}

For commonly used vector norms in $\mathbb{R}^n$, their
dual norms are given by
\[
\| \cdot \|_{1^*} = \| \cdot \|_{\infty}, \quad
\| \cdot \|_{{\infty}^*} = \| \cdot \|_1, \quad
\| \cdot \|_{2^*} = \| \cdot \|_2.
\]

For matrices in $\mathbb{R}^{m \times n}$, we consider
the norm corresponding
to the scalar product $\langle A , B \rangle = \mathrm{trace}
( A^{\mathrm{T}} B )$. Thus, we have $\| A \|_{2^*} =
\| \sigma (A) \|_1$ (see \cite{SS90}), where $\sigma (A)$ is the
vector of the singular values of $A$. Since
$\mathrm{trace} (A^{\mathrm{T}} A) = \| A \|_F^2$, we have
$\| A \|_{F^*} = \| A \|_F$.

For linear operators from $E$ to $G$, $\| \cdot \|_{E,G}$ denotes
the operator norm induced from $\| \cdot \|_E$ and $\| \cdot \|_G$.
Similarly, for linear operators from $G$ to $E$, the norm
induced from the dual norms $\| \cdot \|_{E^*}$ and
$\| \cdot \|_{G^*}$, is denoted by $\| \cdot \|_{G^*,E^*}$.

For the adjoint operators and dual norms, we have the following
result \cite{BG09}:

\begin{lemma} \label{lemmaAdjoint}
\[
\| J \|_{E,G} = \| J^* \|_{G^*,E^*}.
\]
\end{lemma}

In particular, when $G$ has lower dimension than $E$, we can
use the lower dimensional $\| J^* \|_{G^*,E^*}$ instead of the higher
dimensional $\| J \|_{E,G}$.

Now, we consider the product space
$E = E_1 \times ... E_p$, where each Euclidean space
$E_i$ is associated with a scalar product
$\langle \cdot , \cdot \rangle_{E_i}$ and the corresponding norm
$\| \cdot \|_{E_i}$, $1 \le i \le p$. In $E$, we consider
the scalar product defined by
\[
\langle (\mathbf{u}_1,...,\mathbf{u}_p), \
(\mathbf{v}_1,...,\mathbf{v}_p) \rangle =
\langle \mathbf{u}_1 ,\ \mathbf{v}_1 \rangle_{E_1} + \cdots +
\langle \mathbf{u}_p ,\ \mathbf{v}_p \rangle_{E_p}
\]
and the corresponding product norm defined by
\[
\| (\mathbf{u}_1,...,\mathbf{u}_p) \|_v =
v ( \| \mathbf{u}_1 \|_{E_1} ,..., \| \mathbf{u}_p \|_{E_p} ) ,
\]
where $v$ is an absolute norm on $\mathbb{R}^p$, that is,
$v( | \mathbf{x} | ) = v ( \mathbf{x} )$,
for any $\mathbf{x} \in \mathbb{R}^p$, where $| \mathbf{x} | =
[ | x_i | ]$, see \cite{LT85} for details.

Let $v^*$ be the dual norm of $v$ and satisfy the usual inner
product on $\mathbb{R}^p$, we are interested in the dual norm
$\| \cdot \|_{v^*}$ of the product norm $\| \cdot \|_v$
which satisfies the scalar product in $E$. The following result
can be found in \cite{BG09}.

\begin{lemma} \label{lemmaProductNorm}
The dual norm of a product norm can be expressed by
\[
\| (\mathbf{u}_1 ,..., \mathbf{u}_p) \|_{v^*} =
v^* ( \| \mathbf{u}_1 \|_{E_1^*} ,..., \| \mathbf{u}_p \|_{E_p^*} ).
\]
\end{lemma}

With the necessary background in adjoint operators
and dual norms, we apply them to the condition numbers for the weighted
least squares problem. We can think of the Euclidean space $E$ with
norm $\| \cdot \|_E$ as the space of the data in the weighted least
squares problem and $G$ with norm $\| \cdot \|_G$ as the space
of the solution in the weighted least squares problem. Then
the function $g$ in (\ref{eqnMapSol}) is an operator from $E$
to $G$ and the condition number is the measurement of the
sensitivity of $g$ to the perturbation in its input data.

Following \cite{Rice66}, if $g$ is Fr{\'e}chet differentiable in
a neighborhood of $\mathbf{y} \in E$, then the condition number
$\mathcal{K}$ of $g$ at $\mathbf{y}$ is given by
\[
\mathcal{K} = \| g'( \mathbf{y} ) \|_{E,G} =
\max_{\| \mathbf{z} \|_E = 1} \| g'( \mathbf{y} ) \cdot \mathbf{z} \|_G ,
\]
where $\| \cdot \|_{E,G}$ is the operator norm induced from the norms
$\| \cdot \|_E$ and $\| \cdot \|_G$. If $g$ is
nonzero, then we define
\[
\mathcal{K}^{\rm rel} = \mathcal{K} \frac{\| \mathbf{y} \|_E}
{\| g ( \mathbf{y} ) \|_G}
\]
as the relative condition number of $g$ at $\mathbf{y} \in E$.
The above definition shows that $\mathcal{K}$ is dependent of
the norm of the the linear operator $g'( \mathbf{y} )$. Applying
Lemma~\ref{lemmaAdjoint}, we have the following expression of
$\mathcal{K}$ in terms of adjoint operator and dual norm:
\begin{equation} \label{eqnK}
\mathcal{K} = \max_{\| \Delta \mathbf{y} \|_E = 1}
\| g' ( \mathbf{y} ) \cdot \Delta \mathbf{y} \|_G =
\max_{\| \mathbf{z} \|_{G^*} = 1}
\| g' (\mathbf{y})^* \cdot \mathbf{z} \|_{E^*}.
\end{equation}

Now we consider the componentwise measurement on the data space
$E = \mathbb{R}^n$. For any given $\mathbf{y} \in E$,
$E_Y$ denotes the set of all the perturbations
$\Delta \mathbf{y} \in \mathbb{R}^n$ such that
$\Delta y_i = 0$ when $y_i = 0$, $1 \le i \le n$.
Thus in the componentwise perturbation analysis, we use the norm
\[
\| \Delta \mathbf{y} \|_c = \min \{ \omega , \
| \Delta y_i | \le \omega \, |y_i|, \ i=1,...,n \}
\]
to measure the perturbation $\Delta \mathbf{y} \in E_Y$
of $\mathbf{y}$. We call $\| \cdot \|_c$ the componentwise
relative norm. Equivalently,
\begin{equation} \label{eqnCNorm}
\| \Delta \mathbf{y} \|_c = \max \{ |\Delta y_i| / |y_i| \} =
\| (|\Delta y_1| / |y_1| ,..., |\Delta y_n| / |y_n|) \|_{\infty} ,
\end{equation}
where $\Delta \mathbf{y} \in E_Y$.

Next we investigate the dual norm $\| \cdot \|_{c^*}$
of the componentwise norm $\| \cdot \|_c$.
Let the product space $E$
be $\mathbb{R}^n$, each $E_i$ be $\mathbb{R}$,
and the absolute norm $v$ be $\| \cdot \|_{\infty}$.
Setting the norm $\| \Delta y_i \|_{E_i}$ in $E_i$ to
$| \Delta y_i | / | y_i |$ when $y_i \neq 0$, from
Definition~\ref{defDualNorm}, we have the dual norm
\[
\| \Delta y_i \|_{E_i^*} =
\max_{z \neq 0} \frac{| \Delta y_i \cdot z |}{\| z \|_{E_i}} =
\max_{z \neq 0} \frac{| \Delta y_i \cdot z |}{|z| / |y_i|} =
|\Delta y_i| \, |y_i| .
\]
Applying Lemma~\ref{lemmaProductNorm} and (\ref{eqnCNorm}) and recalling
$\| \cdot \|_{\infty^*} = \| \cdot \|_1$, we get the dual norm
\begin{equation} \label{eqnDualC}
\| \Delta \mathbf{y} \|_{c^*} =
\| (\|\Delta y_1\|_{E^*} ,..., \|\Delta y_n\|_{E^*}) \|_{\infty^*} =
\| (|\Delta y_1| \, |y_1| ,..., |\Delta y_n| \, |y_n|) \|_1 .
\end{equation}

Due to the condition $\| \Delta \mathbf{y} \|_E = 1$ in the
condition number $\mathcal{K}$ in (\ref{eqnK}), whether
$\Delta \mathbf{y}$ is in $E_Y$ or not, the expression
of the condition number $\mathcal{K}$ remains valid. Indeed,
if $\Delta \mathbf{y} \not\in E_Y$, that is, $\Delta y_i = 0$
for some $i$, then $\| \Delta \mathbf{y} \|_c = \infty$.
Consequently, such perturbation $\Delta \mathbf{y}$ is
excluded from the calculation of $\mathcal{K}$. Following
(\ref{eqnK}), we have the following lemma on the condition
number in adjoint operator and dual norm.

\begin{lemma} \label{lemmaK}
Using the above notations and the componentwise norm defined
in (\ref{eqnCNorm}), the condition number $\mathcal{K}$ can be
expressed by
\[
\mathcal{K} =
\max_{\| \mathbf{u} \|_{G^*} = 1}
\| (g'(\mathbf{y}))^* \cdot \mathbf{u} \|_{c^*} ,
\]
where $\| \cdot \|_{c^*}$ is given by (\ref{eqnDualC}).
\end{lemma}

Having discussed the norms on the data space, in the next section,
we study the norms on the solution space, which can be either
componentwise or normwise. However, regardless of the norms chosen
in the solution space, we always use the componentwise norm in
the data space.

\section{Condition Numbers for the Weighted Least Squares
Problem} \label{secWLS}

In this section, we present an explicit expression of the
condition number for the weighted least squares problem.
First, we derive an explicit expression of the Fr{\'e}chet
derivative of the mapping $g$ in (\ref{eqnMapSol}), when
$A$ is of full column rank. Let $B \in \mathbb{R}^{m \times n}$,
$\mathbf{c} \in \mathbb{R}^m$ and
$J = g'(A, \mathbf{b})$ be the
derivative, applying the chain rule, we get
\begin{eqnarray}
J(B, \mathbf{c}) &=& g'(A, \mathbf{b}) \cdot
(B, \mathbf{c}) \nonumber \\
&=& L^{\mathrm{T}} (A^{\mathrm{T}}WA)^{-1}
(A^{\mathrm{T}}W (\mathbf{c} - B \mathbf{x}) +
B^{\mathrm{T}}W (\mathbf{b} - A \mathbf{x})) \nonumber \\
&=& L^{\mathrm{T}} ((A^{\mathrm{T}}WA)^{-1} B^{\mathrm{T}} \mathbf{d} -
A_{W,I_n}^{\dagger} B \mathbf{x}) +
L^{\mathrm{T}} A_{W,I_n}^{\dagger} \mathbf{c} ,
\label{eqnJ}
\end{eqnarray}
recalling that $\mathbf{d} = W (\mathbf{b} - A \mathbf{x})$.
Note that $J(B, \mathbf{c})$ is a mapping from
the data space $\mathbb{R}^{m \times n} \times \mathbb{R}^m$ to
$\mathbb{R}^k$.

From the definition of the adjoint operator and
the definition of the scalar
product in the data space $\mathbb{R}^{m \times n} \times \mathbb{R}^m$,
the following lemma gives an explicit
expression of the adjoint operator of the above $J(B, \mathbf{c})$.

\begin{lemma} \label{lemmaDualJ}
The adjoint operator of the Fr{\'e}chet derivative
$J(B, \mathbf{c})$ in (\ref{eqnJ}) is
\[
J^*( \mathbf{u} ) =
(\mathbf{d}\mathbf{u}^{\mathrm{T}} L^{\mathrm{T}}
((A^{\mathrm{T}}WA)^{-1} -
(A_{W,I_n}^{\dagger})^{\mathrm{T}}
L \mathbf{u} \mathbf{x}^{\mathrm{T}} , \
A_{W,I_n}^{\dagger} L \mathbf{u} ) ,
\]
for $\mathbf{u} \in \mathbb{R}^k$. Note that
$J^*(\mathbf{u})$ is a mapping
from $\mathbb{R}^k$ to $\mathbb{R}^{m \times n} \times \mathbb{R}^m$.
\end{lemma}
\noindent
\textbf{Proof}.
Let $J_1(B)$ and $J_2(\mathbf{c})$ be the first and second
terms in the sum (\ref{eqnJ}) respectively.
By the definition of the scalar product in the matrix space,
for any $\mathbf{u} \in \mathbb{R}^k$, we have
\begin{eqnarray*}
\langle \mathbf{u}, J_1(B) \rangle &=&
\mathrm{trace} ( L^{\mathrm{T}} (A^{\mathrm{T}}WA)^{-1}
B^{\mathrm{T}} \mathbf{d} \mathbf{u}^{\mathrm{T}} ) -
\mathrm{trace} ( L^{\mathrm{T}} A_{W,I_n}^{\dagger}
B \mathbf{x} \mathbf{u}^{\mathrm{T}} ) \\
&=& \mathrm{trace} ( \mathbf{d} \mathbf{u}^{\mathrm{T}}
L^{\mathrm{T}} (A^{\mathrm{T}}WA)^{-1} B^{\mathrm{T}} ) -
\mathrm{trace} ( \mathbf{x} \mathbf{u}^{\mathrm{T}}
L^{\mathrm{T}} A_{W,I_n}^{\dagger} B ) \\
&=& \mathrm{trace} ( \mathbf{d} \mathbf{u}^{\mathrm{T}}
L^{\mathrm{T}} (A^{\mathrm{T}}WA)^{-1} B^{\mathrm{T}} ) -
\mathrm{trace} ( (A_{W,I_n}^{\dagger})^\mathrm{T} L
\mathbf{u} \mathbf{x}^{\mathrm{T}} B^{\mathrm{T}} ) \\
&=& \langle \mathbf{d} \mathbf{u}^{\mathrm{T}}
L^{\mathrm{T}} (A^{\mathrm{T}}WA)^{-1} , \ B \rangle -
\langle (A_{W,I_n}^{\dagger})^\mathrm{T} L
\mathbf{u} \mathbf{x}^{\mathrm{T}}, \ B \rangle \\
&=& \langle \mathbf{d} \mathbf{u}^{\mathrm{T}}
L^{\mathrm{T}} (A^{\mathrm{T}}WA)^{-1} -
(A_{W,I_n}^{\dagger})^\mathrm{T} L
\mathbf{u} \mathbf{x}^{\mathrm{T}}, \ B \rangle .
\end{eqnarray*}
For $J_2(\mathbf{c})$, we have
\[
\langle \mathbf{u}, \ J_2(\mathbf{c}) \rangle =
\langle \mathbf{u}, \ L^{\mathrm{T}} A_{W,I_n}^{\dagger}
\mathbf{c} \rangle =
\langle (A_{W,I_n}^{\dagger})^\mathrm{T} L \mathbf{u} ,\
\mathbf{c} \rangle .
\]
Let
\[
J_1^*(\mathbf{u}) =
\mathbf{d} \mathbf{u}^{\mathrm{T}}
L^{\mathrm{T}} (A^{\mathrm{T}}WA)^{-1} -
(A_{W,I_n}^{\dagger})^\mathrm{T} L
\mathbf{u} \mathbf{x}^{\mathrm{T}}
\]
and
\[
J_2^*(\mathbf{u}) =
(A_{W,I_n}^{\dagger})^\mathrm{T} L \mathbf{u} ,
\]
then $\langle J^*(\mathbf{u}),\ (B, \mathbf{c}) \rangle =
\langle (J_1^*(\mathbf{u}), J_2^*(\mathbf{u})) , \
(B, \mathbf{c}) \rangle = \langle \mathbf{u},\ J(B, \mathbf{c}) \rangle$,
which completes the proof. \hfill $\Box$

Having obtained an explicit expression of the adjoint operator
of the Fr{\'e}chet derivative, next we give an explicit expression
of the condition number $\mathcal{K}$ (\ref{eqnK}) in terms
the dual norm in the solution space in the following theorem,
where $\mathrm{vec}(A)$ denotes the vector obtained by stacking
the columns of a matrix $A$, $D_A$ denotes the diagonal matrix
$\mathrm{diag}(\mathrm{vec}(A))$, and $\otimes$ is the Kronecker
product operator \cite{Graham81}.

\begin{theorem} \label{thmK}
The condition number for the full rank weighted least squares
problem can be expressed by
\[
\mathcal{K} = \max_{\| \mathbf{u} \|_{G^*} = 1}
\| [ VD_A \ \
A_{W,I_n}^{\dagger} D_\mathbf{b} ]^{\mathrm{T}} L \mathbf{u} \|_1 =
\| [ VD_A \ \ A_{W,I_n}^{\dagger} D_{\mathbf{b}} ]^{\mathrm{T}} L \|_{G^*,1} ,
\]
where
\begin{equation} \label{eqnV}
V = (A^{\mathrm{T}} W A)^{-1} \otimes \mathbf{d}^{\mathrm{T}} -
\mathbf{x}^{\mathrm{T}} \otimes A_{W,I_n}^{\dagger} .
\end{equation}
\end{theorem}

\noindent
\textbf{Proof}.
Let $\Delta A = [\Delta a_{i,j}]$ and $\Delta \mathbf{b} =
[\Delta b_i]$, then, from (\ref{eqnDualC}), we have
\[
\| ( \Delta A \ \ \Delta \mathbf{b} ) \|_{c^*} =
\sum_{i,j} | \Delta a_{i,j} | \, | a_{i,j} | +
\sum_i | \Delta b_i | \, |b_i|.
\]
Applying Lemma~\ref{lemmaDualJ}, we get
\begin{eqnarray*}
& & \| J^* ( \mathbf{u} ) \|_{c^*} \\
&=& \sum_{j=1}^n \sum_{i=1}^m | a_{i,j} |\,
\left| \left( \mathbf{d} \mathbf{u}^{\mathrm{T}} L^T
(A^{\mathrm{T}} W A)^{-1} - (A_{W,I_n}^{\dagger})^{\mathrm{T}} L \mathbf{u}
\mathbf{x}^{\mathrm{T}} \right)_{i,j} \right| +
\sum_{i=1}^m |b_i| \, \left| \left(
(A_{W,I_n}^{\dagger})^{\mathrm{T}} L \mathbf{u} \right)_i \right| \\
&=& \sum_{j=1}^n \sum_{i=1}^m | a_{i,j} |\,
\left| \left(
d_i ( (A^{\mathrm{T}} W A)^{-1} \mathbf{e}_j)^{\mathrm{T}} -
x_j (A_{W,I_n}^{\dagger} \mathbf{e}_i)^{\mathrm{T}}
\right) L \mathbf{u} \right| +
\sum_{i=1}^m |b_i| \, \left| \left(
A_{W,I_n}^{\dagger} \mathbf{e}_i \right)^{\mathrm{T}}
L \mathbf{u} \right| .
\end{eqnarray*}
It can be verified that $d_i (A^{\mathrm{T}} W A)^{-1} \mathbf{e}_j$
is the $(i+(j-1)\,n)$th column of the $n \times (mn)$ matrix
$(A^{\mathrm{T}} W A)^{-1} \otimes \mathbf{d}^{\mathrm{T}}$
and $x_j A_{W,I_n}^{\dagger} \mathbf{e}_i$ is the $(i+(j-1)\,n)$th
column of the $n \times (mn)$ matrix
$\mathbf{x}^{\mathrm{T}} \otimes A_{W,I_n}^{\dagger}$ in
$V$ (\ref{eqnV}), implying that the above expression
equals
\[
\left\| \left[ \begin{array}{c}
D_A V^{\mathrm{T}} L \mathbf{u} \\
D_{\mathbf{b}} (A_{W,I_n}^{\dagger})^{\mathrm{T}} L \mathbf{u}
\end{array} \right] \right\|_1 =
\left\| [ V D_A \ \ A_{W,I_n}^{\dagger} D_{\mathbf{b}} ]^{\mathrm{T}}
L \mathbf{u} \right\|_1 .
\]
The theorem then follows from Lemma~\ref{lemmaK}. \hfill $\Box$

The following case study discusses some commonly used norms
for the norm in the solution space to obtain some specific
expressions of the condition number $\mathcal{K}$.

\begin{corollary} \label{colKinfty1}
Using the above notations,
when the infinity norm is chosen as the norm in the solution
space $G$, we get
\begin{equation} \label{eqnKinfty1}
\mathcal{K}_{\infty} =
\left\| |L^{\mathrm{T}} V| \mathrm{vec}(|A|) +
|L^{\mathrm{T}} A_{W,I_n}^{\dagger}| \, |\mathbf{b}| \right\|_{\infty} .
\end{equation}
\end{corollary}

\noindent
\textbf{Proof}.  When $\| \cdot \|_G = \| \cdot \|_{\infty}$, the dual norm
$\| \cdot \|_{G^*} = \| \cdot \|_1$. Thus
\begin{eqnarray*}
\mathcal{K}_{\infty} &=&
\left\| [ V D_A \ \ A_{W,I_n}^{\dagger} D_{\mathbf{b}} ]^{\mathrm{T}}
L \right\|_1 \\
&=& \| L^{\mathrm{T}} [ V D_A \ \ A_{W,I_n}^{\dagger} D_{\mathbf{b}} ]
\|_{\infty} \\
&=& \left\| |L^{\mathrm{T}} V| \mathrm{vec}(|A|) +
|L^{\mathrm{T}} A_{W,I_n}^{\dagger}| \, |\mathbf{b}| \right\|_{\infty} .
\qquad \Box
\end{eqnarray*}

The following corollary gives an alternative expression of
$\mathcal{K}_{\infty}$.

\begin{corollary} \label{colKinfty2}
Using the above notations,
when the infinity norm is chosen as the norm in the solution
space $G$, we get
\begin{equation} \label{eqnKinfty2}
\mathcal{K}_{\infty} =
\left\|
\sum_{j=1}^n |L^{\mathrm{T}} (A^{\mathrm{T}} W A)^{-1}
(\mathbf{e}_j \mathbf{d}^{\mathrm{T}} -
x_j A^{\mathrm{T}} W)| \, |A(:,j)| +
|L^{\mathrm{T}} A_{W,I_n}^{\dagger}| \, |\mathbf{b}| \right\|_{\infty} .
\end{equation}
\end{corollary}

\noindent
\textbf{Proof}.
Partitioning
\[
V = [V_1 \ ... \ V_n] ,
\]
where each $V_j$, $1 \le j \le n$, is an $n \times m$ matrix,
we get
\begin{equation} \label{eqnKinfty2a}
\mathcal{K}_{\infty} =
\left\| |L^{\mathrm{T}}V| \mathrm{vec}(|A|) +
|L^{\mathrm{T}} A_{W,I_n}^{\dagger}| \, |\mathbf{b}| \right\|_{\infty} =
\left\| \sum_{j=1}^n |L^{\mathrm{T}}V_j| \, |A(:,j)| +
|L^{\mathrm{T}} A_{W,I_n}^{\dagger}| \, |\mathbf{b}| \right\|_{\infty} .
\end{equation}
Recalling that $d_i (A^{\mathrm{T}}WA)^{-1} \mathbf{e}_j -
x_j A_{W,I_n}^{\dagger} \mathbf{e}_i$ is the $(i+(j-1)n)$th column
of $V$, we have
\[
V_j = (A^{\mathrm{T}}WA)^{-1} (\mathbf{e}_j \mathbf{d}^{\mathrm{T}} -
x_j A^{\mathrm{T}}W) .
\]
The expression (\ref{eqnKinfty2}) is obtained
by substituting $V_j$ in (\ref{eqnKinfty2a}) with the above expression
for $V_j$ and noticing that $A_{W,I_n}^{\dagger} =
(A^{\mathrm{T}}WA)^{-1} A^{\mathrm{T}}W$. \hfill $\Box$

The advantage of the expression (\ref{eqnKinfty2}) over
(\ref{eqnKinfty1}) is the absence of the Kronecker product.
Consequently, its computation requires less memory.
To further reduce the computational cost,
we will propose efficient methods for estimating an upper bound of
$\mathcal{K}_{\infty}$ in Section \ref{sec:est}.

When $W=I_m$, the weighted least squares problem (\ref{eqnWLS})
reduces to the standard least squares problem (\ref{eq:LS})
and the condition number $\mathcal{K}_{\infty}$ in (\ref{eqnKinfty2})
reduces to the condition number
\[
K_{\infty} (L,A,{\bf b})=
\left\|
\sum_{j=1}^n |L^{\mathrm{T}} (A^{\mathrm{T}}  A)^{-1}
(\mathbf{e}_j \mathbf{r}^{\mathrm{T}} -
x_j A^{\mathrm{T}} )| \, |A(:,j)| +
|L^{\mathrm{T}} A^{\dagger}| \, |\mathbf{b}| \right\|_{\infty}
\]
for the standard least squares problem
given by Baboulin and Gratton \cite[(3.4)]{BG09}, noticing
that when $W = I_m$, we have $A_{W,I_n}^{\dagger} = A^\dagger$,
$A^{\mathrm{T}} W A = A^{\mathrm{T}} A$, and
$\mathbf{d} = W (\mathbf{b} - A \mathbf{x}) =
\mathbf{b} - A A^{\dagger} \mathbf{b} = \mathbf{r}$.

Next, we consider the 2-norm and derive an upper bound.

\begin{corollary} \label{colK2}
When the 2-norm is used in the solution space, we have
\begin{equation} \label{eqnK2}
\mathcal{K}_2 \le
\sqrt{k} \, \mathcal{K}_{\infty} .
\end{equation}
\end{corollary}

\noindent
\textbf{Proof}.
When $\| \cdot \|_G = \| \cdot \|_2$, then
$\| \cdot \|_{G^*} = \| \cdot \|_2$. From Theorem~\ref{thmK},
\[
\mathcal{K}_2 = \| [ VD_A\ \
A_{W,I_n}^{\dagger} D_{\mathbf{b}}]^{\mathrm{T}} L \|_{2,1} .
\]
It follows from \cite{GVL96} that for any matrix $B$,
$\| B \|_{2,1} = \max_{\| \mathbf{u} \|_2 = 1} \| B \mathbf{u} \|_1
= \| B \hat{\mathbf{u}} \|_1$, where $\hat{\mathbf{u}} \in \mathbb{R}^k$
is a unit 2-norm vector. Applying $\| \hat{\mathbf{u}} \|_1 \le
\sqrt{k}\, \| \hat{\mathbf{u}} \|_2$, we get
\[
\| B \|_{2,1} = \| B \hat{\mathbf{u}} \|_1 \le
\| B \|_1 \| \hat{\mathbf{u}} \|_1 \le
\sqrt{k}\, \| B \|_1 .
\]
Substituting the above $B$ with $[ VD_A\ \
A_{W,I_n}^{\dagger} D_{\mathbf{b}}]^{\mathrm{T}} L$, we have
\[
\mathcal{K}_2 \le \sqrt{k} \,
\| [ VD_A \ \ A_{W,I_n}^{\dagger} D_{\mathbf{b}}]^{\mathrm{T}} L \|_1 ,
\]
which implies (\ref{eqnK2}). \hfill $\Box$

The above upper bound for $\mathcal{K}_2$ can be obtained by
computing (\ref{eqnKinfty1}) or (\ref{eqnKinfty2}).

So far, we have discussed the various mixed condition
numbers, that is, componentwise norm in the data space and
the infinite norm or 2-norm in the solution space. In the
rest of the section, we study the case of componentwise
condition number,
that is, componentwise norm in the solution space as well.

\begin{corollary} \label{colKc}
Considering the componentwise norm defined by
\begin{equation}\label{eq:comp norm}
\| \mathbf{u} \|_c =
\min \{ \omega , \ | u_i | \le \omega \,
|(L^{\mathrm{T}} \mathbf{x})_i|, \ i=1,...,k \} =
\max \{ |u_i| / |(L^{\mathrm{T}} \mathbf{x})_i|, \ i=1,...,k \} ,
\end{equation}
in the solution space, we have the following three expressions
for the componentwise condition number
\begin{eqnarray*}
\mathcal{K}_c
&=& \| D_{L^{\mathrm{T}} \mathbf{x}}^{-1}
L^{\mathrm{T}} [ VD_A \ \
A_{W,I_n}^{\dagger} D_{\mathbf{b}} ] \|_{\infty} \\
&=& \| |D_{L^{\mathrm{T}} \mathbf{x}}^{-1} |
(|L^{\mathrm{T}} V| \mathrm{vec} (|A|) +
|L^{\mathrm{T}} A_{W,I_n}^{\dagger}|\,|\mathbf{b}| \|_{\infty} \\
&=& \left\|
\sum_{j=1}^n |D_{L^{\mathrm{T}} \mathbf{x}}^{-1}
L^{\mathrm{T}} (A^{\mathrm{T}} W A)^{-1}
(\mathbf{e}_j \mathbf{d}^{\mathrm{T}} -
x_j A^{\mathrm{T}} W)| \, |A(:,j)| +
|D_{L^{\mathrm{T}} \mathbf{x}}^{-1}
L^{\mathrm{T}} A_{W,I_n}^{\dagger}| \, |\mathbf{b}| \right\|_{\infty} .
\end{eqnarray*}
\end{corollary}

\noindent
\textbf{Proof}.
The expressions immediately follow from Theorem~\ref{thmK}
and Corollaries \ref{colKinfty1} and \ref{colKinfty2}. \hfill $\Box$

Similarly to $\mathcal{K}_{\infty}$,
when $W = I_m$, the above condition number $\mathcal{K}_c$
reduces to the standard least squares condition number
\[
K_{\infty} (L,A,{\bf b})= \left\|
\sum_{j=1}^n |D_{L^{\mathrm{T}} \mathbf{x}}^{-1}
L^{\mathrm{T}} (A^{\mathrm{T}}  A)^{-1}
(\mathbf{e}_j \mathbf{r}^{\mathrm{T}} -
x_j A^{\mathrm{T}} )| \, |A(:,j)| +
|D_{L^{\mathrm{T}} \mathbf{x}}^{-1}
L^{\mathrm{T}} A^{\dagger}| \, |\mathbf{b}| \right\|_{\infty}
\]
presented in \cite[page 15]{BG09}.

\section{Condition Number Estimators} \label{sec:est}

In this section, we propose efficient methods for estimating
$\mathcal{K}_\infty $ and $\mathcal{K}_c$,
when integrated into the Paige's method \cite{Paige79a, Paige79b}
for solving the weighted least squares problem.

Firstly, we give upper bounds for $\mathcal{K}_\infty $ and
$\mathcal{K}_c$ in the following theorem. 

\begin{theorem}\label{thm:upp}
Using the notations above, we have the upper bounds
\[
\mathcal{K}_{\infty} \leq \mathcal{K}_{\infty}^u :=
\left\|
L^{\mathrm{T}} (A^{\mathrm{T}} W A)^{-1} D_{ |A|^{\rm T} | \mathbf{d} |}
\right\|_{\infty} +
\left\|
L^{\mathrm{T}} A_{W,I_n}^{\dagger} D_{ |A|\, |\mathbf{x}| }
\right\|_{\infty} +
\left\|
L^{\mathrm{T}} A_{W,I_n}^{\dagger} D_{ |\mathbf{b}|}
\right\|_{\infty}
\]
and
\begin{eqnarray*}
\mathcal{K}_{c} \leq \mathcal{K}_{c}^u &:=&
\left\|
D_{L^{\mathrm{T}} \mathbf{x}}^{-1}
L^{\mathrm{T}} (A^{\mathrm{T}} W A)^{-1}
D_{|A|^{\rm T} |\mathbf{d}|}
\right\|_{\infty} +
\left\|
D_{L^{\mathrm{T}} \mathbf{x}}^{-1} L^{\mathrm{T}} A_{W,I_n}^{\dagger}
D_{|A|\, |\mathbf{x}|}
\right\|_{\infty} \\
& & + \mbox{} \left \|
D_{L^{\mathrm{T}} \mathbf{x}}^{-1}
L^{\mathrm{T}} A_{W,I_n}^{\dagger} D_{|\mathbf{b}|}
\right\|_{\infty}.
\end{eqnarray*}
\end{theorem}

\noindent
\textbf{Proof}.
From the monotonicity property of infinity norm and triangle
inequality, we get
\begin{eqnarray*}
& & \mathcal{K}_{\infty} \\
&\leq& \left\|
\sum_{j=1}^n \left(
|L^{\mathrm{T}} (A^{\mathrm{T}} W A)^{-1}| \mathbf{e}_j \,
|\mathbf{d}^{\mathrm{T}}| \, |A(:,j)| +
| x_j  L^{\mathrm{T}} (A^{\mathrm{T}} W A)^{-1} A^{\mathrm{T}} W)| \,
|A(:,j)| \right) +
|L^{\mathrm{T}} A_{W,I_n}^{\dagger}| \, |\mathbf{b}|
\right\|_{\infty} \\
&\leq& \left\|
\sum_{j=1}^n \left( |L^{\mathrm{T}} (A^{\mathrm{T}} W A)^{-1}| \,
\mathbf{e}_j \, |A(:,j)|^{\rm T}  |\mathbf{d}| +
| x_j  L^{\mathrm{T}} A_{W,I_n}^{\dagger} | \, |A(:,j)| \right)
\right\|_{\infty} +
\left \|
|L^{\mathrm{T}} A_{W,I_n}^{\dagger}| \, |\mathbf{b}|
\right\|_{\infty} \\
&=& \left\|
|L^{\mathrm{T}} (A^{\mathrm{T}} W A)^{-1}| \,
|A|^{\rm T} |\mathbf{d} | \right\|_{\infty} +
\left\| |L^{\mathrm{T}} A_{W,I_n}^{\dagger}| \, |A|\, |\mathbf{x}|
\right\|_{\infty} +
\left\|
|L^{\mathrm{T}} A_{W,I_n}^{\dagger}| \, |\mathbf{b}| \right\|_{\infty} \\
&=& \left\|
L^{\mathrm{T}} (A^{\mathrm{T}} W A)^{-1}
D_{|A|^{\rm T} |\mathbf{d}|} \right\|_{\infty} +
\left\| L^{\mathrm{T}} A_{W,I_n}^{\dagger} D_{|A|\, |\mathbf{x}|}
\right\|_{\infty} +
\left \|
L^{\mathrm{T}} A_{W,I_n}^{\dagger} D_{|\mathbf{b}|}
\right\|_{\infty},
\end{eqnarray*}
where the last equation can be obtained by applying
\[
\|BD_{\mathbf{v}}\|_{\infty} =
\| \,|BD_{\mathbf{v}}|\,\|_{\infty}
= \| \,|B|\,|D_{\mathbf{v}}| \,\|_{\infty} =
\| \,|B| \,|D_{\mathbf{v}}|{\bf 1} \|_{\infty}
= \left\| \,|B| \,|\mathbf{v}| \right\|_{\infty}
\]
where ${\bf 1}=[1,...,1]^{\rm T}$.

The upper bound of $\mathcal{K}_c$ can be derived similarly.
\hfill $\Box$

Our experiments show that the above upper bounds are tight.

The above upper bounds can be computed efficiently when the
weighted least squares problem is solved by the fast
numerically stable method proposed by Paige \cite{Paige79a, Paige79b}.
To see this, we briefly describe the Paige's method.

The weighted least squares problem (\ref{eqnWLS}) arises in
finding the least squares estimate of the vector $\mathbf{x}$
in the linear model $\mathbf{b} = A \mathbf{x} + \mathbf{w}$,
where $A$ is an $m \times n$ matrix and $\mathbf{w}$ is an unknown
noise vector of zero mean and $m \times m$
covariance $Z = W^{-1}$. Usually, the factorization $Z = BB^{\rm T}$
is available. Paige considers the following form equivalent to the
weighted least squares (\ref{eqnWLS}):
\begin{equation} \label{eqnWLSa}
\min_{\mathbf{v}, \mathbf{x}} \| \mathbf{v} \|_2^2 \quad
\hbox{subject to } \mathbf{b} = A \mathbf{x} + B \mathbf{v}.
\end{equation}
By applying the plane rotations, we can get a generalized QR factorization
\cite[(2.1)]{Paige79a} of the data matrix $[\mathbf{b} \ \ A \ \ B]$
of (\ref{eqnWLSa}):
\begin{equation} \label{eqnGQR}
Q^{\rm T} [\mathbf{b} \ \ A \ \ B]
\begin{array}{c}
\begin{array}{ccc}
1 & \, n \, & m
\end{array} \\
\left[ \begin{array}{ccc}
1 & 0   & 0 \\
0 & I_n & 0 \\
0 & 0   & P
\end{array} \right]
\end{array} =
\begin{array}{cc}
  &  \begin{array}{ccccc} \ \ 1\,  &  \,  n\, & m-n-1 & 1\  & n\  \ \end{array} \\
\begin{array}{r} m-n-1 \\ 1 \\ n \end{array} &
\left[ \begin{array}{ccccc}
     0     &    0      & \quad \  L_1 \ \ \quad & 0      & 0 \\
   \eta    &    0      & \mathbf{g}^{\rm T} &     \rho   & 0 \\
\mathbf{z} & R^{\rm T} &      L_{21}        & \mathbf{s} & L_2
\end{array} \right]
\end{array}
\end{equation}
where $Q, P \in \R^{m\times m}$ are orthogonal matrices
and $L_1$, $L_2$, $R^{\rm T}$ are lower triangular and nonsingular,
and $\rho$ is nonzero, assuming $A$ is of full column rank and
$B$ is nonsingular and lower triangular. It is shown that the
weighted least squares solution $\mathbf{x}$ can be obtained by
solving the following nonsingular lower triangular system:
\[
\left[ \begin{array}{cc}
\rho & 0 \\
\mathbf{s} & R^{\rm T}
\end{array} \right] \,
\left[ \begin{array}{c}
\mu \\
\mathbf{x}
\end{array} \right] =
\left[ \begin{array}{c}
\eta \\
\mathbf{z}
\end{array} \right] .
\]
The cost of Paige's algorithm is $O(m^2n/2 + mn^2 -2n^3/3)$
\cite[(4.4)]{Paige79a}.

Letting
\[
S = \left[ \begin{array}{ccc}
     L_{1}         &     0      & 0 \\
\mathbf{g}^{\rm T} &   \rho     & 0 \\
      L_{21}       & \mathbf{s} & L_{2}
\end{array} \right] ,
\]
then
\[
S^{-1} = \left[ \begin{array}{ccc}
               L_1^{-1}                &      0        & 0 \\
-\rho^{-1} \mathbf{g}^{\rm T} L_1^{-1} &   \rho^{-1}   & 0 \\
L_2^{-1} (\rho^{-1} \mathbf{sg}^{\rm T} - L_{21}) L_1^{-1} &
-\rho^{-1} L_2^{-1} \mathbf{s} & L_2^{-1}
\end{array} \right] .
\]
From (\ref{eqnGQR}), we have
\[
A = Q \left[ \begin{array}{c}
0 \\
R^{\rm T}
\end{array} \right] \quad \hbox{and} \quad
B = QSP^{\rm T} .
\]
Consequently, we get the factorizations:
\begin{eqnarray}
W &=& Z^{-1} = QS^{-{\rm T}} S^{-1} Q^{\rm T}, \nonumber \\
A^{\rm T} W A &=& RL_2^{-{\rm T}}L_2^{-1}R^{\rm T}, \quad
(A^{\rm T} W A)^{-1} = R^{-{\rm T}}L_2 L_2^{\rm T}R^{-1}, \label{eq:decomp} \\
A_{W,I_n}^{\dagger} &=&
(A^{\rm T} W A)^{-1} A^{\rm T} W =
= R^{-{\rm T}}
\left[ \rho^{-1} \mathbf{sg}^{\rm T} L_1^{-1} - 
L_{21}L_1^{-1} \ \ \rho^{-1} \mathbf{s} \ \ I_{n} \right]
Q^{\rm T}. \nonumber
\end{eqnarray}
Note that $L_1$, $L_2$,
$R$, and $S$ are triangular matrices.
Then, using the above factorizations,
the three norms in $\mathcal{K}_{\infty}^u$ or $\mathcal{K}_{c}^u$
can be efficiently estimated by the classical condition estimation
method \cite[Chapter 15]{Higham02},
as shown in Algorithm \ref{al:mixed}.

\begin{algorithm} \label{al:mixed}
Estimating $\mathcal{K}_{\infty}^u$.
\begin{algorithmic}
\STATE
Initial vectors $\mathbf{h}_i = k^{-1} {\bf 1} \in \R^k$,
$i=1,2,3$;
\FOR {$p=1,2,\ldots$}
\STATE
Using \eqref{eq:decomp}, calculate
\[
\mathbf{y}_1 = D_{|A|^{\rm T} |\mathbf{d}|}  (A^{\mathrm{T}} W A)^{-1}
L \, \mathbf{h}_1; \ 
\mathbf{y}_2 = D_{|A|\, |\mathbf{x}| } (A_{W,I_n}^{\dagger})^{\rm T}
L \, \mathbf{h}_2; \ 
\mathbf{y}_3 = D_{|\mathbf{b}|} (A_{W,I_n}^{\dagger})^{\rm T}
L \, \mathbf{h}_3;
\]
\STATE
Compute $\mathbf{s}_i ={\rm sign}(\mathbf{y}_i)$, $i=1,2,3$,
where {\rm sign} is the sign function;
\STATE
Using \eqref{eq:decomp}, calculate
\[
\mathbf{z}_1 = L^{\mathrm{T}} (A^{\mathrm{T}} W A)^{-1} \,
D_{|A|^{\rm T} |\mathbf{d}|} \mathbf{s}_1; \ 
\mathbf{z}_2 = L^{\mathrm{T}} A_{W,I_n}^{\dagger}
D_{|A|\, |\mathbf{x}|} \mathbf{s}_2; \
\mathbf{z}_3 = L^{\mathrm{T}} A_{W,I_n}^{\dagger}
D_{|\mathbf{b}|} \mathbf{s}_3;
\]
\IF
{$\|\mathbf{z}_i\|_\infty \leq \mathbf{h}_i^{\rm T} \mathbf{z}_i$}
\STATE
$\gamma_i = \left\| \mathbf{y}_i \right\|_1$, $i=1,2,3$;
\STATE
break
\ENDIF
\STATE
$\mathbf{h}_i = \mathbf{e}_{k_i}$, where $k_i$ is the smallest index
such that $|z_{k_i}|=\| \mathbf{z}_i \|_\infty$;
\ENDFOR
\STATE
Return $\hat{\mathcal{K}}_{\infty}^u = \gamma_1 + \gamma_2 + \gamma_3$.
\end{algorithmic}
\end{algorithm}

Table \ref{tab:al1} lists the major costs in
Algorithm \ref{al:mixed} when integrated into the Paige's method,
where $\mathbf{v}$ is a vector
with conformal dimensions. Let $p_{\max}$ be the total number
of iterations, then the total cost of Algorithm \ref{al:mixed}
is $O(p_{\max}(m^2 + mn + n^2))$. Recalling that the cost of the
Paige's method for solving the weighted least squares problem
is $O(m^2n/2 + mn^2 -2n^3/3)$.

\begin{table}
\caption{
Major operations and their costs of Algorithm \ref{al:mixed}.
} \label{tab:al1}
\centering
\begin{tabular}{cc}
\hline
Operations & flops \\
\hline
$|A^{\rm T}| |\mathbf{d}|$      & $(2m-1)n$ \\
$|A|\, |\mathbf{x}|$            & $(2n-1)m$ \\
$L\mathbf{h}_i$                 & $(2k-1)n$ \\
$Q \mathbf{v}$                  & $(2m-1)m$ \\
$L_{21}^{\rm T} \mathbf{v}$     & $(2n-1)(m-n-1)$ \\
$\mathbf{s}^{\rm T} \mathbf{v}$ & $2n-1$ \\
$R^{-1} \mathbf{v}$ ($R^{-{\rm T}} \mathbf{v}$) & $O(n^2)$ \\
$L_2 \mathbf{v}$ ($L_2^{\rm T} \mathbf{v}$)     & $O(n(n+1))$ \\
$L_1^{-{\rm T}}(\rho^{-1} \mathbf{gs}^{\rm T} -
L_{21}^{\rm T}) \mathbf{v}$ & $O(mn+(m-n-1)^2)$ \\
$(\rho^{-1} \mathbf{sg}^{\rm T} -L_{21})
L_1^{-1} \mathbf{v}$        & $O(mn+(m-n-1)^2)$\\
\hline
\end{tabular}
\end{table}

\section{An Error Analysis of the Augmented System}
\label{secAugmented}

In this section we perform a componentwise perturbation
analysis of the augmented system (\ref{eqnAugmented})
for the weighted least squares problem. Our analysis is
a generalization of the analysis of the standard least
square problem by Arioli et al. \cite{ADR89} and
Bj\"{o}rck \cite{Bjorck91}. 

Let the perturbations $\Delta A \in \mathbb{R}^{m \times n}$ and
$\mathbf{b} \in \mathbb{R}^m$ satisfy
$|\Delta A| \le \epsilon \, |A|$ and
$|\mathbf{b}| \le \epsilon \, |\mathbf{b}|$ for a small $\epsilon$.
Suppose that the perturbed augmented system is
\[
\left[ \begin{array}{cc}
           W^{-1}           & A + \Delta A \\
(A + \Delta A)^{\mathrm{T}} & 0
\end{array} \right] \,
\left[ \begin{array}{c}
\mathbf{d} + \Delta \mathbf{d} \\
\mathbf{x} + \Delta \mathbf{x}
\end{array} \right] =
\left[ \begin{array}{c}
\mathbf{b} + \Delta \mathbf{b} \\
\mathbf{0}
\end{array} \right].
\]
Denoting
\[
G = \left[ \begin{array}{cc}
W^{-1} & A \\
A^{\mathrm{T}} & 0
\end{array} \right], \quad
\mathbf{f} = \left[ \begin{array}{c}
\mathbf{b} \\ \mathbf{0}
\end{array} \right], \quad
\mathbf{z} = \left[ \begin{array}{c}
\mathbf{d} \\ \mathbf{x}
\end{array} \right] ,
\]
and the perturbations
\[
\Delta G = \left[ \begin{array}{cc}
0 & \Delta A \\
(\Delta A)^{\mathrm{T}} & 0
\end{array} \right], \quad
\Delta \mathbf{f} = \left[ \begin{array}{c}
\Delta \mathbf{b} \\ \mathbf{0}
\end{array} \right], \quad
\Delta \mathbf{z} = \left[ \begin{array}{c}
\Delta \mathbf{d} \\ \Delta \mathbf{x}
\end{array} \right].
\]
When $A$ is of full column rank, $G$ is invertible.
It can be verified that
\[
G^{-1} = \left[ \begin{array}{cc}
W - A_{W,I_n}^{\dagger} A^{\mathrm{T}} W &
A_{W,I_n}^{\dagger} \\
(A_{W,I_n}^{\dagger})^{\mathrm{T}} &
- (A^{\mathrm{T}} W A)^{-1}
\end{array} \right] .
\]
We know that if the spectral radius
\begin{equation} \label{eqnSpectral}
\rho \left( |G^{-1}| \, |\Delta G|\right) < 1
\end{equation}
then $I_{m+n} + G^{-1} \Delta G$ is invertible. Clearly, the condition
\begin{equation} \label{eqnEpsilon}
\epsilon <
\rho^{-1} \left( \left[ \begin{array}{cc}
|A_{W,I_n}^{\dagger}|\, |A|^{\mathrm{T}} &
|W - A_{W,I_n}^{\dagger} A^{\mathrm{T}} W|\,|A| \\
|(A^{\mathrm{T}} W A)^{-1}|\,|A|^{\mathrm{T}} &
|(A_{W,I_n}^{\dagger})^{\mathrm{T}}|\,|A|
\end{array} \right] \right) ,
\end{equation}
implies \eqref{eqnSpectral}.
The following results \cite{Skeel79} are necessary
for Theorem \ref{thmAugmented}.

\begin{lemma} \label{lemmaPerturb}
For a linear system $G \mathbf{z} = \mathbf{f}$ and its
perturbed system
\[
(G + \Delta G)(\mathbf{z} + \Delta \mathbf{z}) =
\mathbf{f} + \Delta \mathbf{f} ,
\]
where $\mathbf{z} + \Delta \mathbf{z}$ is the solution
to the perturbed system,
when the perturbations $\Delta G$ and $\Delta \mathbf{f}$
are sufficiently small such that $G+\Delta G$ is invertible,
the perturbation $\Delta \mathbf{z}$
in the solution $\mathbf{z}$ satisfies
\[
\Delta \mathbf{z} = (I + G^{-1} \Delta G)^{-1} G^{-1}
(\Delta \mathbf{f} - \Delta G \mathbf{z}),
\]
which implies
\[
| \Delta \mathbf{z} | \le
\left| (I + G^{-1} \Delta G)^{-1} \right| \, | G^{-1} |
(| \Delta \mathbf{f} | + | \Delta G | \, |\mathbf{z}|).
\]
Furthermore, when the spectral radius
$\rho(| G^{-1} | \, | \Delta G |) < 1$, we have
\begin{eqnarray}
| \Delta \mathbf{z} | &\le&
(I - |G^{-1}|\,|\Delta G|)^{-1} | G^{-1} |
(| \Delta \mathbf{f} | + | \Delta G | \, |\mathbf{z}|) \nonumber \\
&=& (I + O(|G^{-1}|\,|\Delta G|)) | G^{-1} |
(| \Delta \mathbf{f} | + | \Delta G | \, |\mathbf{z}|). \label{eqnDz}
\end{eqnarray}
\end{lemma}

Now we have the bounds for the perturbations in the weighted
least squares solution and residual.

\begin{theorem} \label{thmAugmented}
Using the above notations, for any $\epsilon > 0$ satisfying
the condition (\ref{eqnEpsilon}),
when the componentwise perturbations
$|\Delta A| \le \epsilon \, |A|$ and
$|\Delta \mathbf{b}| \le \epsilon \, |\mathbf{b}|$,
the error in the solution is bounded by
\begin{equation} \label{eqnDx}
\| \Delta \mathbf{x} \|_{\infty} \le \epsilon \left(
\| |(A_{W,I_n}^{\dagger})^{\mathrm{T}}|
(|\mathbf{b}| + |A|\,|\mathbf{x}|) \|_{\infty} +
\|(A^{\mathrm{T}} W A)^{-1}|\,|A|^{\mathrm{T}} |\mathbf{d}| \|_{\infty}
\right) + O(\epsilon^2)
\end{equation}
and error in the weighted residual is bounded by
\begin{equation} \label{eqnDd}
\| \Delta \mathbf{d} \|_{\infty} \le \epsilon \left(
\| |W - A_{W,I_n}^{\dagger} A^{\mathrm{T}} W|
(|\mathbf{b}| + |A|\,|\mathbf{x}|) \|_{\infty} +
\| |(A_{W,I_n}^{\dagger})^{\mathrm{T}}|\,
|A|^{\mathrm{T}} |\mathbf{d}| \|_{\infty}
\right) + O(\epsilon^2) .
\end{equation}
\end{theorem}

\noindent
\textbf{Proof}.
Since the condition (\ref{eqnEpsilon}) implies (\ref{eqnSpectral}),
applying (\ref{eqnDz}) in Lemma~\ref{lemmaPerturb}, we get
\[
\left[ \begin{array}{c}
\Delta \mathbf{d} \\
\Delta \mathbf{x}
\end{array} \right] \le
(I + O(|G^{-1}|\,|\Delta G|)) |G^{-1}|
\left[ \begin{array}{c}
|\Delta \mathbf{b}| + |\Delta A|\,|\mathbf{x}| \\
|\Delta A|^{\mathrm{T}} |A|\,|\mathbf{d}|
\end{array} \right] .
\]
Finally, by using the conditions $|\Delta A| \le \epsilon\, |A|$
and $|\Delta \mathbf{b}| \le \epsilon\, |\mathbf{b}|$ and
the explicit form of $G^{-1}$, the upper bounds (\ref{eqnDx})
and (\ref{eqnDd}) can be obtained. \hfill $\Box$

\section{Numerical Experiments} \label{secExp}

In this section, we present our experimental results to
demonstrate the effectiveness of our condition numbers and
their estimators for the weighted LS problem.
All the numerical experiments were carried
out in \textsc{Matlab} 2015b, with the machine precision $\mu \approx
2.2 \times 10^{-16}$.

Firstly, we adopted the example in Baboulin and Gratton \cite{BG09}
and modified $A$, $W$, and $\mathbf{b}$ as the following:
\begin{eqnarray*}
A &=& \left[ \begin{array}{ccc}
1 & 1 & \epsilon^2 \\
\epsilon & 0 & \epsilon^2 \\
0 & \epsilon & \epsilon^2 \\
\epsilon^2 & \epsilon^2 & 2
\end{array} \right] ,\quad
W = U^{\mathrm{T}} \left[ \begin{array}{cccc}
1 & 0        & 0      & 0 \\
0 & 10\gamma & 0      & 0 \\
0 & 0        & \gamma & 0 \\
0 & 0        & 0      & \gamma / 10
\end{array} \right] U ,\\
{\bf b} &=& {\bf b}_1+10^{-5}\cdot{\bf b}_2,\quad
{\bf b}_1 =\left[ \begin{array}{c} 3 \epsilon \\
\epsilon^2 + \epsilon \\
\epsilon^2 + \epsilon \\
2/ \epsilon + 2 \epsilon^3
\end{array} \right],\quad {\bf b}_2=
\left[ \begin{array}{c}
-\epsilon + \epsilon^4 \\
1 - \epsilon^4 / 2 \\
1 - \epsilon^4 / 2 \\
-\epsilon^2 + \epsilon^3 / 2
\end{array} \right] ,
\end{eqnarray*}
where $\epsilon , \gamma > 0$ and
$U$ is a random orthogonal matrix obtained from the
QR decomposition of a random matrix.
As we can see, $A$  (or $W$) becomes ill-conditioned as $\epsilon$
(or $\gamma$) decreases to zero. The vector $\mathbf{b}$
is constructed so that the solution is imbalanced, that is, its
components range widely, to show the benefit of the componentwise
condition. Note that ${\bf b}_1 \in {\rm  Rang}(A)$ and
${\bf b}_2\in {\rm Ker }(A^{\rm T})$, where
${\rm  Rang}(A)$ and ${\rm Ker }(A^{\rm T})$
denote the range space of $A$ and the null space of
$A^{\rm T}$, respectively. We generated the perturbations:
\begin{equation}\label{eqnPert}
\Delta A= 10^{-8} \cdot E \odot A \hbox{ and }
\Delta {\bf b}=10^{-8} \cdot {\bf f} \odot {\bf b},
\end{equation}
where entries of $E$ and ${\bf f}$ are random variables uniformly
distributed in the interval $ (-1,\, 1)$.  Thus the perturbation
size $\| (\Delta A, \Delta \mathbf{b}) \|_c \approx 10^{-8}$.
For the $L$ matrix in our condition numbers, we chose
\[
L_0=I_3,\
L_1=\left[\begin{matrix} 1 & 0 \\
                          0 & 1 \\
                          0 & 0\end{matrix}\right], \hbox{ and }
L_2=\left[\begin{matrix} 0 & 0 & 1\end{matrix}\right]^{\mathrm{T}}.
\]
That is, corresponding to the above three matrices, the whole
$\mathbf{x}$, the subvector $[x_1 \ x_2]^{\mathrm{T}}$, and
the component $x_3$ are selected respectively.
We called the \textsc{Matlab} built-in function $\tt lscov$ to
compute the solutions $\bfx$ and $\tilde{\mathbf{x}}$ corresponding
to the unperturbed WLS \eqref{eqnWLS} and its perturbed WLS defined by
\[
(A + \Delta A)^{\mathrm{T}} W (A + \Delta A) \tilde{\mathbf{x}} =
(A + \Delta A)^{\mathrm{T}} W (\mathbf{b} + \Delta \mathbf{b}).
\]

We measured the mixed and componentwise relative errors in
$L^{\rm T} \mathbf{x}$ defined by
\[
\mathcal{E}_\infty^{\rm rel} =
\frac{\|L^{\mathrm{T}} \tilde{\mathbf{x}} -
L^{\mathrm{T}} \mathbf{x} \|_\infty}
{\|L^{\mathrm{T}} \mathbf{x} \|_\infty} \hbox{ and }
\mathcal{E}_c^{\rm rel} =
\frac{\|L^{\mathrm{T}} \tilde{\mathbf{x}} -
L^{\mathrm{T}} \mathbf{x} \|_c}
{\|L^{\mathrm{T}} \mathbf{x} \|_c},
\]
where $\|\cdot\|_c$ is the componentwise norm defined
in \eqref{eq:comp norm}. Since the data perturbation size
is about $10^{-8}$, $\mathcal{E}_\infty^{\rm rel} \times
10^{8}$ and $\mathcal{E}_c^{\rm rel} \times 10^{8}$
are respectively indications of the mixed and componentwise
condition numbers for this particular problem.
Specifically, in the table, our condition numbers are
\begin{eqnarray*}
\mathcal{K}_\infty^{\rm rel} &=&
\left\| |L^{\mathrm{T}} V|{\rm vec}(|A|) +
|L^{\mathrm{T}} A_{W,I_n}^{\dagger}|\,|b| \right\|_\infty /
\left\| L^{\mathrm{T}} \mathbf{x} \right\|_\infty,\\
\mathcal{K}_c^{\rm rel} &=&
\left\| |D^{-1}_{L^{\mathrm{T}} \mathbf{x}}
|(|L^{\mathrm{T}} V|{\rm vec}(|A|) +
|L^{\mathrm{T}} A_{W,I_n}^{\dagger}|\,|\mathbf{b}|) \right\|_\infty,
\end{eqnarray*}
where $V$ is defined in \eqref{eqnV}.
Also, in the table,
we define the upper bound $\mathcal{K}_\infty^{u,{\rm rel}} =
\mathcal{K}_\infty^{u} / \| L^{\rm T} \mathbf{x} \|_2$,
recalling that $\mathcal{K}_{\infty}^u$ and $\mathcal{K}_c^u$
are defined in Theorem~\ref{thm:upp}.

\begin{table}
\caption{
Comparison of our condition numbers $\mathcal{K}_\infty^{\rm rel}$
and $\mathcal{K}_c^{\rm rel}$ and their upper bounds
$\mathcal{K}_{\infty}^{u,\rm rel}$
and $\mathcal{K}_c^u$ with their corresponding relative
errors $\mathcal{E}_\infty^{\rm rel}$ and
$\mathcal{E}_c^{\rm rel}$.
} \label{tab:atableold}
\centering
\begin{tabular}{ccccccccc}
\hline
$\epsilon$ & $\gamma$ & $L$ & $\mathcal{E}_\infty^{\rm rel}$ &
$\mathcal{K}_{\infty}^{\rm rel}$ & $\mathcal{K}_\infty^{u,{\rm rel}}$ &
$\mathcal{E}_c^{\rm rel}$ & $\mathcal{K}_c^{\rm rel}$ &
$\mathcal{K}_c^{u}$ \\
 \hline
$10^{-2}$ & $10^0$ & $L_0$ &
4.0597e$-09$ & 2.0000e$+00$ & 2.0000e$+00$ & 1.3431e$-08$ &
4.0806e$+02$ & 4.0813e$+02$ \\
 & & $L_1$ &
1.1355e$-08$ & 2.6350e$+02$ & 3.4499e$+02$ & 1.3431e$-08$ &
4.0806e$+02$ & 4.0813e$+02$ \\
& & $L_2$ &
4.0597e$-09$ & 2.0000e$+00$ & 2.0000e$+00$ & 4.0597e$-09$ &
2.0000e$+00$ &2.0000e$+00$ \\
 \hline
$10^{-2}$ & $10^{-6}$ & $L_0$ &
8.5775e$-09$ & 2.0004e$+00$ & 2.0005e$+00$ & 1.2050e$-06$ &
5.3548e$+02$ & 5.3566e$+02$ \\
& & $L_1$ &
9.5653e$-07$ & 3.3222e$+02$ & 4.2338e$+02$ &
1.2050e$-06$ & 5.3548e$+02$ & 5.3566e$+02$ \\
& & $L_2$ &
8.5775e$-09$ & 2.0004e$+00$ & 2.0004e$+00$ &
8.5775e$-09$ & 2.0004e$+00$ & 2.0004e$+00$ \\
 \hline
$10^{-6}$ & $10^0$ & $L_0$ &
1.2311e$-07$ & 9.1456e$+00$ & 9.1456e$+00$ &
1.9008e$-02$ & 1.4121e$+06$ & 1.4121e$+06$ \\
& & $L_1$ &
1.9008e$-02$ & 9.9851e$+05$ & 1.4121e$+06$ &
1.9008e$-02$ & 1.4121e$+06$ & 1.4121e$+06$ \\
& & $L_2$ &
1.0898e$-09$ & 2.0000e$+00$ & 2.0000e$+00$ &
1.0898e$-09$ & 2.0000e$+00$ & 2.0000e$+00$ \\
 \hline
$10^{-6}$ & $10^{-6}$ & $L_0$ &
5.2538e$-09$ & 4.6298e$+03$ & 4.6298e$+03$ &
8.1423e$-04$ & 1.0280e$+09$ & 1.0280e$+09$ \\
& & $L_1$ &
8.1423e$-04$ & 7.2690e$+08$ & 1.0280e$+09$ &
8.1423e$-04$ & 1.0280e$+09$ & 1.0280e$+09$ \\
& & $L_2$ &
5.2538e$-09$ & 2.0000e$+00$ & 2.0000e$+00$ &
5.2538e$-09$ & 2.0000e$+00$ & 2.0000e$+00$ \\
\hline
\end{tabular}
\end{table}

Table \ref{tab:atableold} compares our condition numbers
$\mathcal{K}_{\infty}^{\rm rel}$ and $\mathcal{K}_c^{\rm rel}$
with their corresponding relative errors $\mathcal{E}_{\infty}^{\rm rel}$
and $\mathcal{E}_c^{\rm rel}$. First, the table shows that our condition
numbers, mixed and componentwise, are consistently close
to the estimates $\mathcal{E}^{\rm rel} \times 10^{8}$. Second, our
componentwise condition numbers show that the third component
of the solution is better conditioned than the first two,
showing the benefit of the componentwise analysis. Third,
in the case when $\epsilon = \delta = 10^{-6}$, our condition
numbers are much larger than their corresponding estimates
$\mathcal{E}_{\infty}^{\rm rel} \times 10^{8}$ and
$\mathcal{E}_c^{\rm rel} \times 10^{8}$.
Our explanation is that
$\mathcal{E}_{\infty}^{\rm rel} \times 10^{8}$ and
$\mathcal{E}_c^{\rm rel} \times 10^{8}$
give estimates of the condition numbers for this particular problem
with this particular perturbation, whereas our condition numbers
are upper bounds for this problem with general perturbation.

Secondly, we experimented on the linear model:
\[
\mathbf{b} = A \mathbf{x} + \mathbf{w},
\]
where $\mathbf{x} \in \mathbb{R}^n$ whose entries are random
variables with standard normal distribution and
$\mathbf{w} \in \mathbb{R}^m$ whose entries $w_i$ are random variables
with normal distribution, mean 0, and predefined variances
$\sigma_i^2$. Thus the weight matrix $W = D_{\bf z}^{-1}$,
where ${\bf z} = [ \sigma_i^2 ]$. In our experiments,
we set $m = 50$ and $n = 10$. The $m \times n$ matrix $A$
was generated by the \textsc{Matlab} built-in function
$\mathtt{sprandn}$ with density $0.5$. The same as before,
the perturbations on $A$ and $\mathbf{b}$ were generated by
(\ref{eqnPert}) and both the unperturbed and perturbed
weighted least squares problems were solved by the
\textsc{Matlab} function $\mathtt{lscov}$.

For the $L$ matrix in our condition numbers, we chose
\[
L_0=I_n,\quad
L_1=\left[ \begin{array}{cc}
1 & 0 \\
0 & 1 \\
\vdots & \vdots \\
0 & 0 \\
0 &0
\end{array} \right]\in \R^{n\times 2},\quad
L_2= \mathbf{e}_n.
\]
Thus, corresponding to the above three matrices, the whole
$\mathbf{x}$, the subvector $[x_{1} \ x_2]^{\mathrm{T}}$, and
the last component $x_n$ are selected respectively.

To investigate the impact of the variances $\sigma_i^2$,
we first set $\mathbf{z} =
[10^{-4}:(4 \times 10^{-4})/(m-1):5 \times 10^{-4}]$,
that is, $\sigma_i^2$ are evenly spaced between $10^{-4}$ and
$5 \times 10^{-4}$.
We generated 1,000 samples of $A$ and $\mathbf{b}$ each. The mean values of
$\mathcal{E}_\infty^{\rm rel}$, $\mathcal{E}_c^{\rm rel}$,
$\mathcal{K}_\infty^{\rm rel}$ and $\mathcal{K}_c^{\rm rel}$,
$\mathcal{K}_\infty^{u,{\rm rel}}$ and $\mathcal{K}_c^{u}$
are displayed in Table~\ref{tblEx1}.
As expected, the condition numbers are moderate, when all the
variances are small, that is, the weight matrix $W$ is well-conditioned,
and the data matrix $A$ is well conditioned since it is a random matrix.

\begin{table}
\caption{
Comparison of our condition numbers $\mathcal{K}_\infty^{\rm rel}$
and $\mathcal{K}_c^{\rm rel}$ and their upper bounds $\mathcal{K}_\infty^{u,\rm rel}$
and $\mathcal{K}_c^{u}$ with the corresponding relative
errors $\mathcal{E}_\infty^{\rm rel}$ and $\mathcal{E}_c^{\rm rel}$,
when the variances $\sigma_i^2$ are evenly spaced between
$10^{-4}$  and $5 \times 10^{-4}$.
} \label{tblEx1}
\begin{center}
\begin{tabular}{ccccccc}
\hline
$L$ & $\mathcal{E}_\infty^{\rm rel}$ & $\mathcal{K}_{\infty}^{\rm rel}$ &
$\mathcal{K}_\infty^{u, {\rm rel}}$ & $\mathcal{E}_c^{\rm rel}$ &
$\mathcal{K}_c^{\rm rel}$ & $\mathcal{K}_c^{u}$  \\
 \hline
$L_0$ & 5.5085e$-09$ & 2.7060e$+00$ & 4.5323e$+00$ &
1.5998e$-07$ & 2.6098e$+02$ & 2.6204e$+02$ \\ 
$L_1$ & 6.8834e$-09$ & 7.5328e$+00$ & 8.4972e$+00$ &
3.9877e$-08$ & 6.2049e$+01$ & 6.2311e$+01$ \\ 
$L_2$ & 1.9161e$-08$ & 2.8946e$+01$ & 2.9075e$+01$ &
1.9161e$-08$ & 2.8946e$+01$ & 2.9075e$+01$ \\ 
 \hline
\end{tabular}
\end{center}
\end{table}

We then widened the range of the variances. Specifically,
$\sigma_i^2$ are evenly spaced between $10^{-4}$ and
$10^2$. Table~\ref{tblEx4} shows the average values of
$\mathcal{E}_\infty^{\rm rel}$, $\mathcal{E}_c^{\rm rel}$,
$\mathcal{K}_\infty^{\rm rel}$ and $\mathcal{K}_c^{\rm rel}$,
$\mathcal{K}_\infty^{u,{\rm rel}}$ and $\mathcal{K}_c^{u}$
over 1,000 samples of $A$ and $\mathbf{b}$ each.
As the range of the variances
is widened, that is, the condition number of the weight matrix $W$
increases, the condition numbers of the weighted least squares
problem increase. However, both Tables~\ref{tblEx1} and \ref{tblEx4},
along with Table~\ref{tab:atableold}, show that the condition
of the weighted least squares problem is more sensitive to
the condition of the data matrix $A$ than that of
the weight matrix $W$.

\begin{table}
\caption{
Comparison of our condition numbers  $\mathcal{K}_\infty^{\rm rel}$
and $\mathcal{K}_c^{\rm rel}$ and their upper bounds $\mathcal{K}_\infty^{u,\rm rel }$
and $\mathcal{K}_c^{u}$ with their corresponding relative
errors $\mathcal{E}_\infty^{\rm rel}$ and $\mathcal{E}_c^{\rm rel}$,
when the variances $\sigma_i^2$ are evenly spaced between $10^{-4}$
and $10^2$.
} \label{tblEx4}
\begin{center}
\begin{tabular}{ccccccc}
\hline
$L$ & $\mathcal{E}_\infty^{\rm rel}$ & $\mathcal{K}_{\infty}^{\rm rel}$ &
$\mathcal{K}_\infty^{u, {\rm rel}}$ & $\mathcal{E}_c^{\rm rel}$ &
$\mathcal{K}_c^{\rm rel}$ & $\mathcal{K}_c^{u}$ \\
 \hline
$L_0$ & 9.2329e$-09$ & 6.4432e$+00$ & 8.8514e$+00$ &
6.0428e$-07$ & 1.0567e$+03$ & 1.2003e$+03$ \\ 
$L_1$ & 1.2621e$-08$ & 1.6575e$+01$ & 1.7419e$+01$ &
1.3424e$-07$ & 2.4916e$+02$ & 2.8492e$+02$ \\ 
$L_2$ & 4.8094e$-08$ & 7.6436e$+01$ & 8.6338e$+01$ &
4.8094e$-08$ & 7.6436e$+01$ & 8.6338e$+01$ \\
 \hline
\end{tabular}
\end{center}
\end{table}

Our experiments show that our upper bounds are consistently very close to
their corresponding condition numbers. In other words, our condition
number estimators are accurate as well as efficient.

\section*{Conclusion} \label{secConclusion}

By applying adjoint operator and dual norm theory, we define
the mixed and componentwise
condition numbers for the linear solution function of the
weighted linear least squares problem. Both the normwise and
componentwise perturbation analyses of the solution
are performed. Moreover,
we present the componentwise perturbation analysis of both
the solution and the residual of the augmented system of the
weighted least squares problem.
We also propose two efficient condition number estimators.
Our numerical experiments show that our condition numbers are
tight and can reveal the condition numbers of individual components
of the solution. Moreover, our condition number estimators are
accurate as well as efficient.

%


\end{document}